\documentclass[]{amsart}
\usepackage{amsmath,amsthm, amsfonts, amssymb, color,graphicx}
\newtheorem{thm}{\sc Theorem}
\newtheorem{cor}{\sc Corollary}
\newtheorem{lem}{\sc Lemma}

\theoremstyle{definition}

\definecolor{wco}{rgb}{0.5,0.2,0.3}

\def\be#1{\begin{equation}\label{#1}}
\def\ee{\end{equation}}

 \def\dist{{\sf dist}}\def\diam{{\rm diam\;}}
\def\C{\mathbb{C}}\def\N{\mathbb{N}}\def\R{\mathbb{R}}\def\Z{\mathbb{Z}}
\def\ds{\displaystyle}
\def\Y{{\mathcal Y}}\def\P{{\mathcal P}}
\def\Q{{\mathcal Q}}\def\W{{\mathcal W}}
\def\ende{{\sf q.e.d.}}\def\f{\mathfrak{f}}
\def\Ende{\ende\medskip}
\def\proof{\noindent{\sc Proof.~}}
\def\remark{\noindent{\sc Remark.~}}
\title{The Yosida Class Is Universal}
\author{Norbert Steinmetz}

\begin{document}
\maketitle

\today
\bigskip{\small\begin{quote}{\sc Abstract.} We discuss
families of meromorphic functions $f_h$ obtained from single
functions $f$ by the re-scaling process
$f_h(z)=h^{-\alpha}f(h+h^{-\beta}z)$ generalising Yosida's process
$f_h(z)=f(h+z)$. The main objective is to obtain information on
the value distribution of the generating functions $f$ themselves.
Among the most prominent generalised Yosida functions are first,
second and fourth Painlev\'e transcendents. The Yosida class
contains all limit functions of generalised Yosida functions--the
Yosida class is universal.\end{quote}

\medskip\noindent{\sc Keywords.} Normal family, Nevanlinna theory, spherical derivative,
Painlev\'e transcendents, elliptic function, Yosida function,
re-scaling,

\medskip\noindent{\sc 2010 MSC.} 30D30, 30D35, 30D45}

\section{Introduction}\label{S1}

\noindent{\sc Yosida functions.} In \cite{KY2}~Yosida introduced
the class $(A)$ of transcendental meromorphic in the complex plane
having bounded spherical derivative
 \be{spherderi}f^\#=\frac{|f'|}{1+|f|^2}.\ee
Then the translates $f_h(z)=f(z+h)$ of $f$ in the class $(A)$ form
a normal family in $\C$, and vice versa; $f$ is called of the {\it
first category}, if no limit function
 \be{limfkt}\f=\lim\limits_{h_n\to\infty}f_{h_n}\ee
is a constant (convergence is always understood with respect to
the spherical metric). It is this additional condition that makes
the class $A_0$ of {\it Yosida functions} so fascinating. The
elementary functions (like $e^z$, $\tan z$ etc) have bounded
spherical derivatives, but are not Yosida functions. On the other
hand $A_0$ contains the elliptic functions. A thorough
investigation of the class $A_0$ was performed by
Favorov~\cite{JuF}, with emphasis on the distribution of zeros and
poles.

\medskip\noindent{\sc Painlev\'e transcendents.} The {\it first
Painlev\'e transcendents} are the solutions to Painlev\'e's first
differential equation $w''=z+6w^2;$ they are meromorphic in $\C$
(see \cite{NSt1}) and satisfy $w^\#=O(|z|^{\frac34}).$ More
precisely, if $\Q$ denotes the set of (non-zero) zeros of $w$ and
$\Q_\epsilon=\bigcup_{q\in\Q}\{z:|z-q|<\epsilon|q|^{-\frac14}\}$,
then $w^\#(z)=O(|z|^{-\frac14})$ holds outside $\Q_\epsilon$,
while $w^\#(q)\asymp|q|^{\frac34}$. The family $(w_h)_{|h|\ge 1}$
with $w_h(z)=h^{-\frac12}w(h+h^{-\frac14}z)$ is normal in $\C$,
and for ``most'' solutions the limit functions
$\lim\limits_{h_n\to\infty}w_{h_n}$ are non-constant. There exist,
however, solutions with large zero- and pole-free regions, and in
that case one has constant solutions $\not\equiv 0,\infty$,
see~\cite{NSt2}. We note, however, that in many applications it is
only required that $f_{h_n}\not\to 0,\infty$.

\medskip\noindent{\sc Definition.} The class $\Y_{\alpha,\beta}$
($\alpha\in\R$, $\beta>-1$) consists of all in $\C$ transcendental
meromorphic functions $f$ in $\C$, such that the family
$(f_h)_{|h|\ge 1}$ of functions
 \be{familie}f_h(z)=h^{-\alpha}f(h+h^{-\beta}z)\ee
(for any or just one determination of $h^{-\alpha}$ and
$h^{-\beta}$) form a normal family in $\C$, and no limit function
(\ref{limfkt}) is constant. Functions in $\Y_{\alpha,\beta}$ are
called {\it generalised Yosida functions.} We also define the
classes $\Y_{\alpha,-1}$ by postulating normality of the family of
functions $f_h(z)=h^{-\alpha}f(h+hz)$ only in
$\C\setminus\{-1\}$, and postpone the analysis of this class to
the last section.

\medskip\noindent{\sc Remarks and Examples.}

\begin{itemize}
 \item  Some of the
results proved in this paper are not new. This, in particular,
concerns theorems in the classes $\W^0_1=\Y_{0,-1}$ (see the last
section) and $A_0=\Y_{0,0}$. The proofs in this paper cover all
parameters $\alpha\in\R$ and $\beta>-1$.
 \item The re-scaling process in the definition of
$\Y_{\alpha,\beta}$ is motivated by and formally related to the
Pang-Zalcman process \cite{ XP1, XP2, LZ1, LZ2}. As far as I know,
particular classes $\Y_{\alpha,\beta}$ with $\alpha\ne 0$ occurred
for the first time, although implicitly, in the paper \cite{NSt2}
on the Painlev\'e transcendents. Of course, this kind of
re-scaling is also not new. It goes back at least to Valiron, but
was even used in Painlev\'e's so-called $\alpha$-method.

 \item It will turn out that $\Y_{\alpha,\beta}$ is contained in the
 class $\W_{2+|\alpha|+\beta}$
 discussed by Gavrilov~\cite{VIG0}: $f\in \W_p$ $(p\ge 1)$ if and only if
 $\sup_{\C}|z|^{2-p}f^\#(z)<\infty$, see also Makhmutov~\cite{SAM1}; $\W_2$ is Yosida's class $(A)$. The class
$\W_p^{(0)}$, also discussed by
 Gavrilov, coincides with $\Y_{0,p-2}$, while the same is true for
 $A_0$ and $\Y_{0,0}$.
 \item $f\in \Y_{\alpha,\beta}$ implies $1/f\in
 \Y_{-\alpha,\beta}$, and $\tilde f(z)=z^af(z^b)$ ($a\in\Z, b\in\N$)
 belongs to $\Y_{a+b\alpha,b+b\beta-1}$; $z^b$ and $z^a$ may be replaced by
 a polynomial $p$ and a rational function $r$, respectively, with $\deg p=b$ and
 $r(z)\sim cz^a$ as $z\to\infty$ ($c\ne 0)$. We mention two simple corollaries:
  \begin{itemize}\item If $\alpha=-a/b$ is rational, then $\tilde f\in\Y_{0,\beta+b\beta-1}.$
  \item If $-1<\beta<0$ and $b$ is sufficiently large, then $b+b\beta-1\ge 0$.\end{itemize}
  It would thus suffice to deal with the cases $\beta=-1$ and
  $\beta\ge 0$, respectively.
 \item To every $n\in\{2,3,4,6\}$ there exists a meromorphic function $f$
 such that $f(z^n)$ is an elliptic function (see Mues~\cite{EM}).
 Thus $f\in\Y_{0,1-1/n}$, and $\tilde f(z)=z^af(z^b)$
 belongs to $\Y_{a,b/n-1}$ ($a\in\Z,$ $b\in\N$).
 \item $f'\in\Y_{\alpha,\beta}$ implies
 $f\in\Y_{\alpha-\beta,\beta}$ for at least one primitive.
 \item ``Most'' of the first, second and fourth Painlev\'e transcendents belong to
 $\Y_{\frac 12,\frac 14}$,  $\Y_{\frac12,\frac12}$ and $\Y_{1,1}$, respectively
 (for ``some'' solutions the second condition is violated, namely those having
large zero- and pole-free regions).
 \begin{itemize}\item[$\circ$] Any first Painlev\'e
 transcendent has a primitive $W$ which also is a first integral: $w'^2=2zw+4w^3-2W$;
 in ``most''cases $W\in\Y_{\frac14,\frac14}$, although $w'^2, zw, w^3\in\Y_{\frac32,\frac14}.$
 \item[$\circ$] The second Painlev\'e equation $w''=a+zw+2w^3$ has a first
 integral $W$: $w'^2=2aw+zw^2+w^4-W$ with $W'=w^2$; since $w^2\in\Y_{1,\frac12}$ (in ``most''
 cases),  $W\in\Y_{\frac12,\frac12}$ follows, although $w'^2, zw, w^4\in\Y_{2,\frac12}.$
 \item[$\circ$] Painlev\'e's fourth equation $2ww''={w'}^2+3w^4+8zw^3+4(z^2-a)w^2+2b$
 also has a first integral $W$:
 ${w'}^2=w^4+4zw^3+4(z^2-a)w^2-2b-4wW$ with $W'=w^2+2zw$ and $W\in\Y_{1,1}$, again only in ``most'' cases.
  \end{itemize}\end{itemize}

\section{Simple Properties}\label{S2}

\begin{thm}\label{THM1} Every $f\in\Y_{\alpha,\beta}$ satisfies
$f^\#(z)=O(|z|^{|\alpha|+\beta})$.\end{thm}

\proof We may assume $\alpha\ge 0$, otherwise would replace $f$ by
$1/f$, noting that $f^\#=(1/f)^\#$ and $1/f\in\Y_{-\alpha,\beta}$.
For $|h|\ge 1$ we have
 \begin{equation}\label{formelthm1}f_h^\#(0)=|h|^{-\alpha-\beta}f^\#(h)\frac{1+|f(h)|^2}{1+|h|^{-2\alpha}|f(h)|^2}
 \ge |h|^{-\alpha-\beta}f^\#(h),\end{equation}
while the left hand side is bounded by Marty's
Criterion.\quad\Ende

\noindent{\sc Remarks.}

\begin{itemize}
 \item The bound $|z|^{|\alpha|+\beta}$ is sharp
 (not only for the Painlev\'e transcendents).
 \item It is  obvious that every limit function
 $\f=\lim\limits_{n\to\infty}f_{h_n}$ belongs to $\W_2$.
 More precisely, $\f^\#$ is bounded by $m_f=\sup\limits_{z\in\C,
 |h|>1}f^\#_h(z)$: if $z_0$ is not a pole of $\f$, then we have
 also $f'_{h_n}\to \f'$ close to $z_0$, hence
 $\f^\#(z_0)=\lim\limits_{n\to\infty}f^\#_{h_n}(z_0)\le m_f.$ At a
 pole of $\f$ we will consider $1/\f$ instead of $\f$ (more in
 Theorem \ref{THM8}).
 \item The limit functions of the Painlev\'e
 families $(w_h)$ are elliptic functions.
 \end{itemize}\medskip

Yosida~\cite{KY2} has shown that given $f\in A_0$ and $\epsilon>0$
there exists some $\delta>0$, such that
$\ds\int_{|z-h|<\epsilon}\!\!\!\!\!\!\!\!\!\!\!\!f^\#(z)\,d(x,y)>\delta$
holds for every $h\in\C$. The analog for $\Y_{\alpha,\beta}$ is
Theorem \ref{THM2} below. For $\beta$ fixed, $|h|>1$ and
$\epsilon>0$ we set
 \be{Delta}\Delta_\epsilon(h)=\{z:|z-h|<\epsilon|h|^{-\beta}\}.\ee

\begin{thm}\label{THM2} For every $f\in\Y_{\alpha,\beta}$ and $\epsilon>0$
we have
 $$\inf_{|h|>1}|h|^{2|\alpha|}\int_{\Delta_\epsilon(h)}\!\!\!\!\!\!\!\!\!\!f^\#(z)^2\,d(x,y)>0
 \quad{\rm and}\quad\inf_{|h|>1}\sup_{z\in\Delta_\epsilon(h)}
f^\#(z)|z|^{|\alpha|-\beta}>0.$$
\end{thm}

\remark The second inequality was proved by Gavrilov~\cite{VIG1} for the class $\Y_{0,-1}$
(which he denoted $\W_1^0$).

\proof The integral in question is
 $\ds I=\int_{|w|<\epsilon}\!\!\!\!\!\!\!\!\!f^\#(h+h^{-\beta}w)^2|h|^{-2\beta}\,d(u,v).$
From
 $$\begin{array}{rcl}
 f^\#(\zeta)|h|^{-\beta}&=&\ds|h|^\alpha
f_h^\#(w)\frac{1+|h|^{-2\alpha}|f(\zeta)|^2}
 {1+|f(\zeta)|^2}\cr&\ge&\ds\min\{1,|h|^{-2\alpha}\} |h|^\alpha
 f_h^\#(w)=|h|^{-|\alpha|}f_h^\#(w)\end{array}$$
($\zeta=h+h^{-\beta}w$, $|h|\ge 1$)  follows $\ds
|h|^{2|\alpha|}I\ge\int_{|w|<\epsilon}\!\!\!\!\!\!\!\!f_h^\#(w)^2\,d(u,v),$
and by definition of $\Y_{\alpha,\beta}$ the right hand side has a
positive infimum with respect to $h$.\quad\Ende

\begin{thm}\label{THM3}Let $f$ be meromorphic in $\C$. Then in order that $f\in\Y_{0,\beta}$
it is necessary and sufficient that
 $$f^\#(z)=O(|z|^{\beta})\quad{\rm and}\quad
 \liminf\limits_{|h|\to\infty}\sup\limits_{{z\in\Delta_\epsilon}(h)}
 f^\#(z)|z|^{-\beta}>0$$
for some $[all]~\epsilon>0$.
\end{thm}

\proof We just have to prove sufficiency. The first condition
ensures that $(f_h)$ is a normal family in $\C$, and the second
guarantees that the limit functions are non-constant:
$\sup\limits_{|z|<\epsilon}\f^\#(z)>0$.\quad\Ende

\medskip\noindent{\sc Definition.} Given
$f\in\Y_{\alpha,\beta}$ we denote by $\P$ and $\Q$ the set of
non-zero poles and zeros of $f$, respectively (if any), and set
(for the definition of $\Delta_\epsilon$ see (\ref{Delta}))
 \be{PQ}\P_\epsilon=\textstyle \bigcup_{p\in\P}\Delta_\epsilon(p)\quad{\rm and}\quad
 \Q_\epsilon=\bigcup_{q\in\Q}\Delta_\epsilon(q).\ee

\begin{thm}\label{THM4} For $f\in\Y_{\alpha,\beta}$ we have
 $$\inf\limits_{q\in\Q}\dist(q,\P)|q|^{\beta}>0\quad{\rm and}\quad
 \inf\limits_{p\in\P}\dist(p,\Q)|p|^{\beta}>0.$$\end{thm}

\proof Take any sequence $(q_n)$ of zeros such that
$\dist(q_n,\P)|q_n|^{\beta}\to\inf\limits_{q\in\Q}\dist(q,\P)|q|^{\beta}$
and $f_{q_n}\to \f\not\equiv const,$ locally uniformly in $\C$.
Then $\f(0)=0$ implies $|\f(z)|<1$ on some disc $|z|<\delta$,
hence $\liminf\limits_{n\to\infty}\dist(q_n,\P)|q_n|^{\beta}\ge
\delta$ by Hurwitz' theorem, this showing that
$\inf\limits_{q\in\Q}\dist(q,\P)|q|^{\beta}\ge \delta>0$.
Concerning the second assertion we just note that
$1/f\in\Y_{-\alpha,\beta}$, so that the notions ``pole'' and
``zero'' may be interchanged.\quad\Ende

\remark We will say that the zeros and poles of $f$ are
$\beta$-{\it separated.} From now on it will be tacitly assumed
that $\Q_\epsilon\cap\P_\epsilon=\emptyset.$

\begin{thm}\label{THM5}Every $f\in\Y_{\alpha,\beta}$ satisfies
 $$\begin{array}{rrcll}{\rm (i)}&|f(z)|&=&O(|z|^\alpha)&(z\notin\P_\epsilon);\cr
 {\rm (ii)}&|1/f(z)|&=&O(|z|^{-\alpha})&(z\notin\Q_\epsilon);\cr
 {\rm (iii)}&|f(z)|&\asymp&|z|^\alpha&(z\notin\P_\epsilon\cup\Q_\epsilon);\cr
 {\rm (iv)}&|f'(z)/f(z)|&=&O(|z|^{\beta})& (z\notin\P_\epsilon\cup\Q_\epsilon);\cr
 {\rm (v)}&f^\#(z)&=&O(|z|^{\beta-|\alpha|})&(z\notin\P_\epsilon\cup\Q_\epsilon).
 \end{array}
 $$\end{thm}

\proof Let $(h_n)$ be any sequence outside $\P_\epsilon$, such
that $f_{h_n}$ tends to $\f\not\equiv const$, locally uniformly in
$\C$, and $|f(h_n)||h_n|^{-\alpha}$ tends to
$M_\epsilon=\sup\limits_{z\not\in\P_\epsilon}|f(z)||z|^{-\alpha}$.
Then $M_\epsilon=|\f(0)|$ is finite. The second assertion follows
from $1/f\in\Y_{-\alpha,\beta}$, and together we obtain
$|f(z)|\asymp|z|^\alpha$ $(z\notin\P_\epsilon\cup\Q_\epsilon).$
(iv) follows from
$\ds\frac{h_n^{-\beta}f'(h_n)}{f(h_n)}\to\frac{\f'(0)}{\f(0)}\ne\infty,$
and from (iii) and (iv) follows $\ds
f^\#(z)=\frac{|f'(z)|}{|f(z)|} \frac1{|f(z)|+\frac1{|f(z)|}}=
O(|z|^{\beta}|z|^{-|\alpha|})$, hence (v).\quad\Ende

\remark The symbol $\asymp$ has proved very useful: $\phi(z)\asymp\psi(z)$ in some
real or complex region  means $|\phi(z)|=O(|\psi(z)|)$ and
$|\psi(z)|=O(|\phi(z)|)$.

\begin{cor}\label{COR1} Every function $f\in\Y_{\alpha,\beta}$ has
infinitely many zeros and poles.\end{cor}

\proof If $f$ had only finitely many poles, then $f$ were rational
as follows from $f(z)=O(|z|^\alpha)$ outside $\P_\epsilon$, hence
in $|z|>R$.\quad\Ende

\begin{thm}\label{THM6} For $f\in\Y_{\alpha,\beta}$ and $\tilde f\in\Y_{\tilde\alpha,\beta}$
with sets of poles and zeros $\Q$ and $\tilde\Q$, and $\P$ and
$\tilde\P$, respectively, the product $f\tilde f$ belongs to
$\Y_{\alpha+\tilde\alpha,\beta}$ if  $\Q\cup\tilde\Q$ and
$\P\cup\tilde\P$ are $\beta$-separated. In particular, $f^m$
belongs to $\Y_{m\alpha,\beta}$.\end{thm}

\proof The hypotheses ensure that zeros [poles] of $f$ cannot
collide with poles [zeros] of $\tilde f$, hence $f_{h_n}\tilde
f_{h_n}\to\f\tilde\f$.\quad\Ende

By Theorem \ref{THM1} the zeros and poles of
$f\in\Y_{\alpha,\beta}$ are $\beta$-separated. On the other hand,
zeros and poles are {\it equally $\beta$-distributed} in the
following sense:

\begin{thm}\label{THM7}Given $f\in\Y_{\alpha,\beta}$ there exist positive numbers
$\epsilon_0$, $\eta_0$, and  $M$, such that

 \begin{itemize}\item[(i)] every disc $\Delta_{\eta_0}(z_0)$ contains at least one
 zero and one pole;
 \item[(ii)]every disc $\Delta_{\epsilon_0}(z_0)$ contains at most $M$ zeros (counted by
 multiplicities) and no pole, or at most $M$ poles of $f$ and no zeros.
 \end{itemize}
 In particular, the zeros and poles of $f$ have bounded
 multiplicities.\end{thm}

\proof Suppose there exist sequences $h_n\to\infty$ and
$\eta_n\to\infty$, such that $\Delta_{\eta_n}(h_n)$ contains no
poles (the same for zeros), while $f_{h_n}\to \f\not\equiv const$,
locally uniformly in $\C$. Then by Hurwitz' Theorem, $\f$ is
finite in every euclidian disc $|z|<\eta_n$, hence is an entire
function, this contradicting Corollary \ref{COR1}. Similarly, if
we assume that the pair $(\epsilon_0,M)$ does not exist, then
there exist sequences $\epsilon_n\to 0$ and $h_n\to\infty$, such
that $f$ has at least $n$ zeros (say) in
$\Delta_{\epsilon_n}(h_n)$, while $f_{h_n}$ tends to some
non-constant function $\f$. By Hurwitz' theorem, $\f$ has a zero
at the origin of order $\ge n$ for every $n$, which is absurd.
Thus there exists $\epsilon>0$, such that the number of zeros in
$\Delta_{\epsilon_0}(z_0)$ is  bounded, uniformly with respect to
$z_0$. Diminishing $\epsilon_0$, if necessary, it we may achieve
by Theorem \ref{THM1} that none of the discs
$\Delta_{\epsilon_0}(z_0)$ contains a pole.\quad\Ende

\remark It is not hard to prove that there also exists some
$\theta_0>0$, such that $f$ assumes {\it every} value in {\it
every} disc $\Delta_{\theta_0}(z_0)$.

\begin{thm}\label{THM8}{\sc [The Yosida class is universal]} For $f\in\Y_{\alpha,\beta}$
$(\beta>-1)$ the limit functions
$\f=\lim\limits_{h_n\to\infty}f_{h_n}$ belong to the Yosida class
$A_0=\Y_{0,0}$.\end{thm}

\proof First of all $\f$ has bounded spherical derivative, hence
$\f\in \W_2$ and the family $(\f_h)_{h\in\C}$ of translations
$\f_h(z)=\f(z+h)$ is normal in $\C$. Also the corresponding sets
$\P$ and $\Q$ are $0$-separated (euclidian distance between $\P$
and $\Q$ is positive), and equally $0$-distributed: there exist
positive numbers $\epsilon_0$, $\eta_0$ and $M$, such that every
disc $|z-h|<\eta_0$ contains at least one zero and one pole, while
every disc $|z-h|<\epsilon_0$ contains at most $M$ poles [zeros],
and no zeros [poles]. If $(h_n)$ is any sequence tending to
$\infty$, then the disc $|z-h_n|<\eta_0$ contains at least one
zero $q_n$ and one pole $q_n$. Since $|p_n-q_n|\ge 2\epsilon_0$,
 all limit functions of $(\f_{h_n})$ also
have at least one zero and one pole in $|z|\le 2\eta_0$, and
therefore are non-constant.\quad\Ende

\section{Value Distribution}\label{S3}

In this section we are concerned with the value distribution of
functions $f\in\Y_{\alpha,\beta}$. For the definition of the
Nevanlinna functions $T(r,f)$, $m(r,f)$ and $N(r,f)$, and for
basic results in Nevanlinna Theory the reader is referred to
Hayman~\cite{WH} and Nevanlinna~\cite{RN}. From the
Ahlfors-Shimizu formula
$$T(r,f)=\ds\frac
 1\pi\int_0^r\int_{|z|<t}f^\#(z)^2\,d(x,y)\,dt+O(1)$$
and Theorem \ref{THM1} follows
$T(r,f)=O(r^{2(|\alpha|+\beta+1)})$, hence $f\in\Y_{\alpha,\beta}$
has order of growth
 $$\ds\varrho(f)=\limsup\limits_{r\to\infty}
 \frac{\log T(r,f)}{\log r}$$
at most $2(|\alpha|+\beta+1)$. Replacing $f$ by $\tilde
f(z)=z^af(z^b)$ with $\tilde f\in\Y_{a+b\alpha,b+b\beta-1}$ yields
$\varrho(f)=\varrho(\tilde f)/b\le
2(|a+b\alpha|+b\beta+b)/b=2(|a/b+\alpha|+\beta+1),$ and since
$\inf\limits_{a\in\Z,b\in\N}|a/b+\alpha|=0$, we obtain in any
case:

\begin{thm}\label{THM9} Every $f\in\Y_{\alpha,\beta}$ has order of growth $\varrho(f)\le
2\beta+2$.\end{thm}

\remark For the first Painlev\'e transcendents the first estimate
yields $\varrho(w)\le\frac72$, while the order is
$\varrho(w)=\frac52=2(\frac14+1)$(\footnote{There is a misprint in
\cite{NSt2}: ``$T(r,f)=2T(r,w)+O(\log r)$'' for
$f(z)=z^{-1}w(z^2)$, of course, has to be replaced by
``$T(r,f)=T(r^2,w)+O(\log r)$''.}). Similarly we have the (sharp)
estimates $\varrho(w)\le 3=2(\frac12+1)$ and $\varrho(w)\le
4=2(1+1)$ for the second and fourth Painlev\'e transcendents,
respectively (see \cite{NSt2}).

\begin{thm}\label{THM10}Every $f\in\Y_{\alpha,\beta}$
$(\beta>-1)$ has  $\asymp r^{2\beta +2}$ zeros and poles in  $|z|<r:$
 \be{anzahlfkt}n(r,0)\asymp r^{2\beta+2}\quad{\rm and}\quad n(r,\infty)\asymp
 r^{2\beta+2}.\ee
In particular, $f$ has order of growth
$\varrho(f)=2\beta+2$.
\end{thm}

\remark We remind the reader that $\phi(r)\asymp\psi(r)$ means
$\phi(r)=O(\psi(r))$ and $\psi(r)=O(\phi(r))$ as
$r\to\infty$.\medskip

\proof With every pole $p$ in $|p|<r$ we associate the disc
$\Delta_{\epsilon_0}(p)$; by Theorem \ref{THM7} it contains at
most $M$ poles. Starting with $p_1$ ($|p_1|<r$), let $p_2$
($|p_2|<r$) be any of the poles not contained in
$\Delta_{\epsilon_0}(p_1)$, $p_3$ ($|p_3|<r$) not contained in
$\Delta_{\epsilon_0}(p_1)\cup\Delta_{\epsilon_0}(p_2)$, and so
forth; we may arrange that
$|p_\nu|^{-\beta}\ge|p_{\nu+1}|^{-\beta}$ holds. Then obviously
$n(r,\infty)=O(\phi(r))$, where $\phi(r)$ counts how many {\it
mutually disjoint} discs $\Delta_{\epsilon_0/2}(p)$ may be placed
in a large euclidian disc $|z|<r+\epsilon_0r^{-\beta}$. The
geometric answer is $\phi(r)=O(r^{2\beta+2})$, if $\beta\ge 0$,
and $\phi(r)\le\phi(r/2)+O(r^{2\beta+2})$ if $-1<\beta<0$, which
also implies $\phi(r)=O(r^{2\beta+2})$ (consider the radii
$r=2^k$). Thus
 $$n(r,\infty)=O(r^{2\beta +2})$$
holds in any case. To prove the converse, we note that for $r$
sufficiently large the annulus $||z|-r|<\eta_0r^{-\beta}$ contains
at least $c'r^{\beta+1}$ mutually disjoint discs of radius
$\eta_0r^{-\beta}$, hence also at least $c'r^{\beta+1}$ poles.
Again we have to distinguish the cases (i)~$\beta\ge 0$ and
(ii)~$-1<\beta<0$. Starting with $r_1$ sufficiently large we
define in case (i) $r_k=r_{k-1}+2\eta_0r_{k-1}^{-\beta}$, while in
case (ii) $r_k$ denotes the unique solution to the equation
$r_k=r_{k-1}+2\eta_0r_k^{-\beta}$ ($k=2,3,\ldots$); note that
$r\mapsto r-2\eta_0r^{-\beta}$ is increasing on
$r^{-\beta-1}<1/(2|\beta|\eta_0)$ if $-1<\beta<0$. Then the annuli
$\big||z|-r_k\big|<\eta_0r_k^{-\beta}$ are mutually disjoint, and
each contains at least $c'r_k^{\beta+1}$ poles of $f$. We claim
 $$\nu_k=n(r_k,\infty)\ge 2cr_k^{2\beta+2},$$
provided $c$ is sufficiently small, this implying $n(r,\infty)\ge
cr^{2\beta+2}$ for $r$ sufficiently large (note that
$r_k\to\infty$). Assuming $\nu_{k-1}\ge 2cr_{k-1}^{2\beta+2}$ to
be true, we obtain
 $$\begin{array}{rcl}
 \nu_k&\ge& \nu_{k-1}+c'r_{k}^{\beta+1}\ge
 2cr_{k-1}^{2\beta+2}+c'r_{k}^{\beta+1}\cr
 &\ge&
 2cr_{k}^{2\beta+2}-2c(r_k-r_{k-1})(2\beta+2)r_k^{2\beta+1}+c'r_{k}^{\beta+1}\end{array}$$
by the Mean Value Theorem. In case (ii) we have
$r_k-r_{k-1}=2\eta_0r_{k}^{-\beta}$, while in case (i)
$r_k-r_{k-1}=2\eta_0r_{k-1}^{-\beta}\le 3\eta_0r_{k}^{-\beta}$
holds (assuming $r_1$ sufficiently large). We thus obtain
 $$\nu_k\ge 2cr_{k}^{2\beta+2}+r_k^{\beta+1}[c'-2c3\eta_0(2\beta+2)]=2cr_{k}^{2\beta+2}$$
if $c$ is chosen to satisfy $c'-2c3\eta_0(2\beta+2)=0$. Finally
from $r_k=O(r_{k-1})$ follows
 $$r^{2\beta+2}=O(n(r,\infty))$$
in all cases $\beta>-1$. The assertion about the order of growth
now follows from $\varrho(f)\le 2\beta+2$ on one hand, and
$T(r,f)\ge N(r,f)\asymp r^{2\beta+2}$ on the other.\quad\Ende

From the proof we obtain:

\begin{cor}\label{COR2}For  $\beta>-1$ and $c>\eta_0$, every annulus $\big||z|-r|\big|<cr^{-\beta}$
contains $\asymp r^{\beta+1}$ zeros $[poles]$ of
$f\in\Y_{\alpha,\beta}$.\end{cor}

\begin{thm}\label{THM11} For every $f\in\Y_{\alpha,\beta}$ holds
 \be{schmfkt}m(r,f)+m(r,1/f)=\frac 1{2\pi}\int_0^{2\pi}
 \big|\log|f(re^{i\theta})|\big|\,d\theta=O(\log r),\ee
and, in particular,
 \be{charfkt}T(r,f)\sim N(r,f)\sim N(r,1/f)\asymp r^{2\beta+2}.\ee\end{thm}

\remark For the class $\Y_{0,0}$ the first assumption was proved
by Favorov~\cite{JuF}, even with $O(1)$ instead of $O(\log
r)$.\medskip

\proof Let $C_r$ denote the circle $|z|=r$. For
$0<\epsilon<\epsilon_0$ fixed, the contribution of
$C_r\setminus(\Q_\epsilon\cup\P_\epsilon)$ to the integral is
$O(\log r)$ by Theorem \ref{THM4} (it says, among others, that
$|f(z)|\asymp |z|^\alpha$). Let $K$ be a component of
$\Q_\epsilon$ [or $\P_\epsilon$] that intersects the circle $C_r$.
If $K$ contains the zeros $q_\mu$ $(1\le\mu\le m\le M)$ [or the
poles $p_\nu$ $(1\le\nu\le n\le M)$, but not zeros and poles
simultaneously], then
$\Phi(z)=f(z)\prod\limits_{\mu=1}^m(z-q_\mu)^{-1}$ is zero- and
pole- free in $K$ and satisfies $|\Phi(z)|\asymp
r^{\alpha}r^{-m\beta}$ on $\partial K$ by Theorem~\ref{THM4}, and
also in $K$ by the Maximum-Minimum Principle. Thus the
contribution of $C_r\cap K$ to the integral is
 $$\sum_{\mu=1}^m\int_{I_r}\big|\log
 |re^{i\theta}-q_\mu|\big|\,d\theta +O(|I_r|\log r),$$
where $I_r=\{\theta\in [0,2\pi):re^{i\theta}\in K\}$ and $|I_r|$
is its linear measure. From Lemma \ref{LEM2} at the end of section
\ref{S5} follows $|I_r|=O(r^{-\beta-1})$, hence
 $$\int_{I_r}\big|\log |re^{i\theta}-q_\mu|\big|\,d\theta=
 O\Big(\int_0^{r^{-\beta-1}}\!\!\!\!\!\!\!\!\!\!|\log (r\theta)|\,d\theta\Big)
 =O(r^{-\beta-1}\log r).$$
The assertion follows from the fact, that by virtue of Corollary
\ref{COR2} there are at most $O(r^{\beta+1})$ components $K$
intersecting $C_r$.\quad\Ende

\begin{thm}\label{THM12}For every $f\in\Y_{\alpha,\beta}$ and $c\in\C$ we have
$\ds m\Big(r,\frac 1{f-c}\Big)=O(\log r).$
\end{thm}

\remark Yosida~\cite{KY2} proved $m\big(r,1/(f-c)\big)=O(r)$ for
$f\in A_0$.

\medskip\proof We just note that Theorem \ref{THM11} also holds for
$f-c$ instead of $f$. For $\alpha\ge 0$ we have
$f-c\in\Y_{\alpha,\beta}$, while $1/f\in\Y_{-\alpha,\beta}$ if
$\alpha<0$ and
 $$\big|\log|f-c|\big|\le \big|\log|c|\big|+\big|\log|f|\big|+\big|\log|1/f-
 1/c|\big|.\quad\ende$$

\begin{thm}\label{THM13}The $c$-points $(c\ne 0)$ of $f\in\Y_{\alpha,\beta}$
are $\beta$-close to the zeros, and $\beta$-separated from the
poles if $\alpha>0$, and vice versa if $\alpha<0:$
 $$\lim_{\zeta\to\infty, f(\zeta)=c}|\zeta|^\beta\dist(\zeta,\Q)=0\quad{\rm and}
 \quad\inf_{f(\zeta)=c}|\zeta|^\beta\dist(\zeta,\P)>0.$$
For $\alpha=0$ and any pair $(a,b)$ the sets of $a$- and
$b$-points are $\beta$-separated.\end{thm}

\proof The first assertion ($\alpha>0$) follows from
$f-c\in\Y_{\alpha,\beta}$ and Theorem \ref{THM4}. If $(\zeta_n)$
denotes any sequence of $c$-points such that $f_{\zeta_n}\to
\f\not\equiv const$, then we have also
$\zeta_n^{-\alpha}(f(\zeta_n+\zeta_n^{-\beta}z)-c)\to \f(z)$ and
$\f(0)=0$. From Hurwitz' Theorem then follows
$|\zeta_n|^{\beta}\dist(\zeta_n,\Q)\to 0$ $(n\to\infty)$. Finally,
since $\Y_{0,\beta}$ is M\"obius invariant, every pair $(a,b)$ can
play the role of $(0,\infty)$.\quad\Ende

\section{Derivatives}\label{S4}

The derivative of $f_h$ is
$f_h'(z)=h^{-\alpha-\beta}f'(h+h^{-\beta}z)$, and since the limit
functions of the family $(f_h)$ are non-rational, one might expect
that $f'\in\Y_{\alpha+\beta,\beta}.$ Now a trivial necessary
condition for $\phi_n\to \phi\not\equiv const$, locally uniformly
in some domain $D$, is that the $a$-points and $b$-points of
$\phi_n$ are locally uniformly $0$-separated (separated with
respect to euclidian metric in any compact subset of $D$). In
general, $\phi_n\to \phi$ does not imply $\phi_n'\to \phi'$ if
$\phi_n$ has poles, in other word, there is no Weierstrass
Convergence Theorem for meromorphic functions (while the converse
is true: $\phi_n'\to \psi$ implies that $\psi$ has a primitive
$\phi$, and $\phi_n\to \phi+const$). The obstacle that prevents
$\phi_n'$ from converging to $\phi'$ is the existence of colliding
poles of $\phi_n$ and/or of zeros of $\phi_n'$ colliding with
poles.

\begin{lem}\label{LEM1}Suppose that $\phi_n$ converges to $\phi$,
locally uniformly in $|z|<r$, and $\phi$ has a pole of order $m$
at $z=0$. Then for $\phi_n'\to \phi'$, locally uniformly in some
neighbourhood of $z=0$, each of the following conditions is
necessary and sufficient: there exist $\rho>0$ and $n_0$, such
that for $n\ge n_0$
\begin{itemize}\item[(i)] $\phi_n$ has only one pole (of order $m$) in $|z|<\rho$;
\item[(ii)] $\phi_n'$ has no zeros in $|z|<\rho$.\end{itemize}
\end{lem}

\proof Since $\phi'$ has a pole of order $m+1$ at $z=0$, and no
other pole and also no zero in $|z|<2\rho$, it is necessary for
$\phi_n'\to\phi'$, uniformly in some neighbourhood of $z=0$, that
$\phi_n'$ $(n\ge n_0)$ has $m+1$ poles (counted with
multiplicities) and no zero in $|z|<\rho$, say. Since every pole
of $\phi_n$ of order $\ell$ is a pole of order $\ell+1$ of
$\phi_n'$, this means that $\phi_n$ has only one pole in
$|z|<\rho$. Conversely, if $\phi_n$ has only one pole $b_n$ (of
order $m$) with $b_n\to 0$, then we have
$\ds\phi_n(z)=\frac{\psi_n(z)}{(z-b_n)^m},$ $\ds\psi_n\to \psi,$
$\ds\psi(0)\ne 0,$ $\ds\psi_n'\to \psi',$
$z\psi'(z)-m\psi(z)\big|_{z=0}\ne 0$, and
 $$\phi_n'(z)=\frac{(z-b_n)\psi'_n(z)-m\psi_n(z)}{(z-b_n)^{m+1}}\to
 \frac{z\psi'(z)-m\psi(z)}{z^{m+1}}=\phi'(z),$$
uniformly in some neighbourhood of $z=0$. It remains to show that
(ii) implies (i). If $\phi_n$ has $p>1$ different poles in
$|z|<\rho$ of total multiplicity $m$, then by the Riemann-Hurwitz
formula $\phi$ has $m-1$ critical points close to $z=0$, only
$m-p$ of them arising from multiple poles. Thus $\phi_n'$ has
$p-1$ zeros close to $z=0$.\quad\Ende

\remark In any case the sequence $\phi_n'$ tends to $\phi'$,
locally uniformly in $0<|z|<\rho$. If (i) or (ii) is violated,
then some of the poles of $\phi'_n$ collide with zeros of
$\phi_n'$, and in the limit multiplicities disappear as do the
zeros of $\phi_n'$. If $\phi_n=1/\P_n$, $\P_n$ a polynomial of
degree $m$, the equivalence of (i) and (ii) follows from the
Gau{\ss}-Lucas Theorem.

\begin{thm}\label{THM14}In
order that for $f\in\Y_{\alpha,\beta}$ the derivative  $f'$
belongs to $\Y_{\alpha+\beta,\beta}$, each of the following
conditions is necessary and sufficient:
\begin{itemize}
 \item[(i)] $\ds\inf_{p\in\P}|p|^{-\beta}\dist(p,\P\setminus\{p\})>0$;
 \item[(ii)] $\ds\inf_{f'(c)=0}|c|^{-\beta}\dist(c,\P)>0$.
\end{itemize}
\end{thm}

\begin{cor}If the poles of $f\in A_0$ are $0$-separated from each other,
then every derivative of $f$ also belongs to $A_0$.\end{cor}

\noindent{\sc Example.}  We construct $f\in\Y_{0,0}$ such that
$f'\not\in\Y_{0,0}$:
$\ds\phi(z)=\ds\sum\limits_{k=1}^\infty\frac{1}{(z-k^2)^2-k^{-2}}$
is meromorphic in $\C$. If $|z-k^2|\ge k/2$ holds for every $k$,
then $\sum\limits_{k=2}^\infty Mk^{-2}$ is a convergent majorant,
hence $f(z)=o(1)$ as $z\to\infty$ outside $\bigcup_{k\ge
1}\{z:|z-k^2|<k/2\}$, while in case $|z-\ell^2|<\ell/2$ for some
$\ell$ we have $|z-k^2|\ge k/2$ for $k\ne \ell$ and
$f(z)=\ds\frac{1}{(z-\ell^2)^2-\ell^{-2}}+o(1)$ as $z\to\infty$ by
the same reason. Thus the limit functions
$\lim\limits_{h_n\to\infty}\phi_{h_n}$ are either constants or
else have the form $(z-z_0)^{-2}$. Then for $f_0\in\Y_{0,0}$ we
have also $f=f_0+\phi\in\Y_{0,0}$, but $f'\not\in\Y_{0,0}$.

\bigskip From $f'\in\Y_{\alpha+\beta,\beta}$ would follow
$m(r,1/f')=O(\log r)$. This, however, is true anyway and provides
a new proof of Theorem \ref{THM12}.

\begin{thm}\label{THM15}Every $f\in\Y_{\alpha,\beta}$ satisfies
$m(r,1/f')=O(\log r).$\end{thm}

\remark This was proved by Yosida \cite{KY2} for $f\in A_0$ with
$O(\log r)$ replaced by $O(r)$.

\medskip\proof Taking into account that $m(r,1/f)=O(\log r)$ and $m(r,f)=O(\log r)$, hence
also $m(r,f')\le m(r,f)+O(\log r)=O(\log r)$ holds,
\be{yos}m(r,1/f')=-\frac 1{2\pi}\int_0^{2\pi}\log f^\#(re
^{i\theta})\,d\theta+O(\log
 r)\ee
follows.(\footnote{More generally, Yosida~\cite{KY2} proved that
$\ds 2T(r,f)-N_1(r)=\ds-\frac 1{2\pi}\int_0^{2\pi}\log
 f^\#(re^{i\theta})\,d\theta+O(1)$ holds,
where $N_1(r)$ ``counts'' the critical points of $f$. The
following proof is straight forward:
 $$\begin{array}{rcl}\ds\frac 1{2\pi}\int_0^{2\pi}\log f^\#(re ^{i\theta})\,d\theta&=&
 m(r,f')-m(r,1/f')-2m(r,f)+O(1)\cr &=&-[N(r,f)+\overline N(r,f)]+N(r,1/f')-2T(r,f)+2N(r,f)+O(1)\cr
 &=&[N(r,f)-\overline N(r,f)]+N(r,1/f')-2T(r,f)+O(1)\cr
 &=&N_1(r)-2T(r,f)+O(1).\cr\end{array}$$})
We claim that the right hand side of (\ref{yos}) is $O(\log r)$.
The lower estimate follows from Theorem \ref{THM1}: $-\log
f^\#(z)\ge-(|\alpha|+\beta)\log|z|+O(1).$ It remains to prove that
 \begin{equation}\label{absch}-\frac 1{2\pi}\int_0^{2\pi}\log[r^{|\alpha|-\beta}
 f^\#(re^{i\theta})]\,d\theta\le C\end{equation}
holds. To this end we divide $[0,2\pi]$ into $\asymp r^{\beta+1}$
intervals of length $\asymp r^{-\beta-1}$. If (\ref{absch}) is not
true, then there exists a sequence $r_n\to\infty$ and intervals
$I_n$ of length $\asymp r_n^{-\beta-1}$, such that
 $$J_n=-r_n^{\beta+1}\int_{I_n}\log[r_n^{|\alpha|-\beta}
 f^\#(r_ne^{i\theta})]\,d\theta\to\infty.$$
We may assume that $I_n=[-r_n^{-\beta-1},r_n^{-\beta-1}]$ and
$f_{r_n}\to \f\not\equiv const$. From
 $$f^\#_{r_n}(z)=r_n^{-\alpha-\beta}f^\#(\zeta)
 \frac{1+|f(\zeta)|^2}{1+r_n^{-2\alpha}|f(\zeta)|^2}
 \le r_n^{-\alpha-\beta}f^\#(\zeta)\max\{1,r_n^{2\alpha}\}=r_n^{|\alpha|-\beta}f^\#(\zeta)$$
($\zeta=r_ne^{i\theta}=r_n+r_n^{-\beta}z$,
$r_nd\theta=|d\zeta|=r_n^{-\beta}|dz|$) then follows
 $$-\log[r_n^{|\alpha|-\beta}f^\#(\zeta)]\le-\log f^\#_{r_n}(z)\quad{\rm
 and}\quad\limsup_{n\to\infty}J_n\le -\int_{[-i,i]}\log
 \f^\#(z)\,|dz|.\quad\ende$$

\section{Series  And Product Developments}\label{S5}

Since the series $\sum\limits_{p\in\P}|p|^{-s-1}$ and
$\sum\limits_{q\in\Q}|q|^{-s-1}$ diverge if
$s=2\beta+2=\varrho(f)$, and converge if $s>2\beta+2$, the
canonical products $\prod\limits_{q\in\Q}E\big(\frac
zq,[2\beta]+2\big)$ and $\prod\limits_{p\in\P}E\big(\frac
zp,[2\beta]+2\big)$ converge absolutely and locally uniformly;
$E(u,g)=(1-u)e^{u+u^2/2+\cdots+u^g/g}$ denotes the Weierstrass
prime factor of {\it genus} $g$. Hence any $f\in\Y_{\alpha,\beta}$
has the

\medskip\noindent{\sc Hadamard Product Representation}
 $$f(z)=z^se^{S(z)}\frac{\prod_{q\in\Q}E\big(\frac zq,m\big)}
 {\prod_{p\in\P}E\big(\frac zp,m\big)}$$
($s\in\Z$ and $S$ a polynomial with $\deg S\le m=2+[2\beta]$), and
differentiation yields the

\medskip\noindent{\sc Mittag-Leffler Expansion}
 $$\frac{f'(z)}{f(z)}=\frac
sz+S'(z)+\sum_{q\in\Q}\frac{z^m}{(z-q)q^m}-
\sum_{p\in\P}\frac{z^m}{(z-p)p^m}.$$

If we do not insist in {\it absolute} convergence, then much more
can be said.

\begin{thm}\label{THM16} Suppose $f\in\Y_{\alpha,\beta}$ and $f(0)\ne 0,\infty$. Then
\be{summe}\frac{f'(z)}{f(z)}=T_{m-1}(z)+\lim_{r\to\infty}\Big[\sum_{|q|<r}\frac{z^m}{(z-q)q^m}-
 \sum_{|p|<r}\frac{z^m}{(z-p)p^m}\Big]\ee
holds, locally uniformly in $\C\setminus(\P\cup\Q)$; $m$ is any
integer $>\beta$, and $T_{m-1}$ is the {$(m-1)$-th} Taylor
polynomial for $f'/f$ at $z=0$. Each zero $q$ and  pole $p$ in the
sum occurs according to its multiplicity.\end{thm}

\proof The following technique is well-known. Let $\Phi$ be
meromorphic in the plane having simple poles $\xi$ with residues
$\rho(\xi)$, and assume that $\Phi(0)\ne 0,\infty$ and
$|\Phi(z)|=O(|z|^\beta)$ holds on the circles $|z|=r_k\to\infty$.
Then \be{resisatz}I_k(z)=\frac{1}{2\pi i}\int_{|\zeta|=r_k}
 \frac{\Phi(\zeta)z^m}{(\zeta-z)\zeta^m}\,d\zeta=O(r_k^{\beta-m})\to
 0\quad(k\to\infty),\ee
provided $m>\beta$. On the other hand, the Residue Theorem yields
 $$I_k(z)=\Phi(z)+\sum_{|\xi|<r_k}\frac{\rho(\xi)z^m}{(\xi-z)\xi^m}-T_{m-1}(z),$$
with $T_{m-1}$ the $(m-1)$-th Taylor polynomial of $\Phi$ at
$z=0$, hence
 \be{PHISUM}\Phi(z)=T_{m-1}(z)+\lim_{k\to\infty}\sum_{|\xi|<r_k}\frac{\rho(\xi)z^m}{(z-\xi)\xi^m}.\ee
This applies to $\Phi=f'/f$ with poles $p$ and $q$, if $r_k$ can
be chosen to lie outside $\P_\epsilon\cup\Q_\epsilon$. If,
however, $|z|=r_k$ intersects some connected component $C$ of $\P$
and/or $\Q$, we may by virtue of Lemma \ref{LEM2} (see the end of
this section) replace the intersection $C\cap \{z:|z|=r_k\}$ by
one or more subarcs of $\partial C$ of total length
$O(r_k^{-\beta})$. This way we obtain the Jordan curve $\Gamma_k$;
it is contained in the annulus $A_k:\big||z|-r_k\big|<\epsilon
r_k^{-\beta}$ and since there are at most $r_k^{1+\beta}$ such
components, the length of $\Gamma_k$ is $O(r_k)$. To get rid of
$\Gamma_k$ and even $r_k$ we just remark that for $|z|<R$ and
$r_k\to\infty$ we have
 \be{rest}\sum_{q\in A_k\cap\Q}\Big|\frac{z^m}{(z-q)q^m}\Big|+
 \sum_{p\in
 A_k\cap\P}\Big|\frac{z^m}{(z-p)p^m}\Big|=O(r_k^{-m-1}r_k^{\beta+1})\to
 0.\quad\ende\ee

Noting that $\ds\frac{z^m}{(z-\xi)\xi^m}=\frac{1}{z-\xi}+
\sum_{j=0}^{m-1}\frac{z^j}{\xi^{j+1}}= \frac d{dz}\log E\Big(\frac
 z\xi,m\Big),$
we obtain:

\begin{thm}\label{THM17} Every $f\in\Y_{\alpha,\beta}$ may be
written as
 \be{produkt}f(z)=z^se^{S(z)}\lim_{r\to\infty}\frac{\prod_{|q|<r}E\big(\frac z{q},m\big)}
 {\prod_{|p|<r}E\big(\frac z{p},m\big)},\ee
where $m$ is any integer $>\beta$, $s\in\Z$, and $S$ is a
polynomial with $\deg S\le m$. Each zero $q$ and pole $p$ in the
products occurs according to its multiplicity.\end{thm}

\remark There are, of course, also more or less complicated
modifications of Theorem \ref{THM16} if $f$ has multiple poles.
Since $m>\beta$ is arbitrary, the limits
 $$\lim_{r\to\infty}\Big[\sum_{|p|<r }p^{-\mu}-\sum_{|q|<r}q^{-\mu}\Big]$$
exist for any integer $\mu>m>\beta$.\medskip

For $\beta$ an integer, the term in brackets, the sums in
(\ref{rest}), and also the sequence of functions $I_k(z)$ in
(\ref{resisatz}) remain uniformly bounded if we choose
$\mu=m=\beta$, which means that in (\ref{summe}) and
(\ref{produkt}) we may replace $m$ by $\beta$ if we simultaneously
replace $r\to\infty$ by $r_k\to\infty$ for some suitably chosen
sequence $(r_k)$.

\medskip\noindent{\sc Theorem \ref{THM17}a.} {\it If $\beta>-1$ is an integer,
then every $f\in\Y_{\alpha,\beta}$ may be written as
$$f(z)=z^se^{S(z)}\lim_{k\to\infty}\frac{\prod_{|q|<r_k}E\big(\frac z{q},\beta\big)}
 {\prod_{|p|<r_k}E\big(\frac z{p},\beta\big)},$$
for some suitably chosen sequence $r_k\to\infty$; $s$ is an
integer and $S$ is a polynomial with $\deg S\le\beta+1$. In
particular, for $f$ in the Yosida class $A_0$ and also in
$\Y_{\alpha,0}$ this means
$$f(z)=z^se^{az+b}\lim\limits_{k\to\infty}\frac{\prod_{|q|<r_k}\big(1-\frac
z{q}\big)}
 {\prod_{|p|<r_k}\big(1-\frac z{p}\big)}$$
and
$$\frac{f'(z)}{f(z)}=\ds a+\frac
sz+\lim_{k\to\infty}\Big[\sum_{|q|<r_k}\frac{1}{z-q}-
 \sum_{|p|<r_k}\frac{1}{z-p}\Big].$$}

\remark The method of proof of Theorem \ref{THM16} applies also to
$f$ itself. If $f$ has only simple poles $p$ with residues
$\rho(p)$ and if $f(0)\ne \infty$, then
 \be{hilf}f(z)=T(z)+\lim_{k\to\infty}\sum_{p\in
 D_k}\frac{\rho(p)z^m}{(z-p)p^m}\ee
holds, locally uniformly in $\C\setminus\P$; $m$ is any integer
$>\alpha$, and $T$ is the $(m-1)$-th Taylor polynomial of $f$ at
$z=0$. To get rid of $D_k$ we need information about the residues.
Let $C_k$ be any component of $\P_\epsilon$ that intersects the
circle $|z|=r_k$. We may assume that $C_k$ contains the poles
poles $p_k^{(\nu)}$ ($1\le\nu\le n\le M_f$, $p_k=p_k^{(1)}$), and
also that $f_{p_k}\to\f\not\equiv const$. Then the contribution of
the poles in the sequence $(C_k)$ to $\f$ is
 $$\lim_{k\to\infty}\sum_{\nu=1}^{n}
 \frac{\rho(p_k^{(\nu)})p_k^{\beta-\alpha}}{z-(p_k^{(\nu)}-p_k)p_k^\beta}
 =P(z)\prod_{\nu=1}^n(z-a_\nu)^{-1};$$
$P$ is a polynomial of degree $<n$, and the numbers
$a_\nu=\lim\limits_{k\to\infty}(p_k^{(\nu)}-p_k)p_k^\beta$ are not
necessarily distinct, since $\f$ may have multiple poles. If
$a_\kappa=\cdots=a_\lambda$ and $\ne a_\mu$ else, then
$\big|\sum\limits_{\nu=\kappa}^\lambda
\rho(p_k^{(\nu)})\big|=O(|p_k|^{\alpha-\beta})$ holds. Since there
are at most $O(r_k^{\beta+1})$ components $C_k$, the contribution
of the annulus $A_k$ to the sum in (\ref{hilf}) is
$O(r_k^{\alpha-m})\to 0$ as $k\to\infty,$ and again we obtain
\be{fsum}f(z)=T(z)+\lim_{r\to\infty}\sum_{|p|<r}\frac{\rho(p)z^m}{(z-p)p^m}.\ee

In the particular case $f\in A_0=\Y_{0,0}$ with simple poles only,
(\ref{fsum}) holds with $m=\alpha=0$ and $r=r_k$, hence
 $$f(z)=a+\lim_{k\to\infty}\sum_{|p|<r_k}\frac{\rho(p)}{z-p}.$$

We finish this section by proving a technical lemma as follows:

\begin{lem}\label{LEM2}Let $C$ be any domain
that consists of $n$ discs
$\Delta_\epsilon(h_\nu)$ and intersects $|z|=r$. Then for $\epsilon$ sufficiently small
and $r$ sufficiently large, $C$ has diameter and boundary curve
length $\le K_n\epsilon r^{-\beta}$; the constant $K_n$ only
depends on $n$.
\end{lem}

\proof We will prove by induction that there exists some $c>0$,
such that for $\epsilon$ sufficiently small any domain
$C_k=\bigcup_{\kappa=1}^k\Delta_\epsilon(h_\kappa)$ $(1\le k\le
n)$ that intersects $|z|=r$ is contained in the annulus
$A_k:\big||z|-r\big|<2ck\epsilon r$ with $\diam C_k\le 2ck\epsilon
r^{-\beta}$. This is obviously true if $k=1$. Assuming the
assertion to be true for $k$ discs, we consider the domain
$C_{k+1}=C_k\cup\Delta_\epsilon(h)$ satisfying $\diam C_{k+1}\le
2ck\epsilon r^{-\beta}+2\epsilon|h|^{-\beta}$. From
$\Delta_\epsilon(h)\cap A_k\ne\emptyset$ then follows
$|h|^{-\beta}<cr^{-\beta}$, $C_{k+1}\subset A_{k+1}$ and $\diam
C_{k+1}\le 2c(k+1)\epsilon r^{-\beta}$. The limitations imposed on
$\epsilon$ and $c$ are $(1\pm 2cn\epsilon)^{-\beta}\le c$ and
$2cn\epsilon<1$. Thus the diameter of $C$ and the length of the
boundary curve of $C$ is $O(\epsilon r^{-\beta})$.\quad\Ende

\section{The Case $\beta=-1$}\label{S6}

The limit functions of the family of functions
$f_h(z)=h^{-\alpha}f(h+hz)$ have an essential singularity at
$z=-1$, since zeros and poles accumulate there. Hence we postulate
normality only in $\C\setminus\{-1\}$ to define the class
$\Y_{\alpha,-1}$. Apart from this it is not hard to verify that
Theorems \ref{THM1} [$f^\#(z)=O(|z|^{|\alpha|+\beta})$],
\ref{THM2}, \ref{THM3}, \ref{THM4} [$\beta$-separation of $\P$ and
$\Q$], \ref{THM5}, \ref{THM6}, \ref{THM7} [$\varrho(f)=2\beta+2$],
\ref{THM9}, \ref{THM11} [$m(r,f)=O(\log r)$], \ref{THM12},
\ref{THM13}, \ref{THM14},  and \ref{THM15}, as well as Corollary
\ref{COR1} remain true also if $\beta=-1$. Beyond the fact
$\varrho(f)=0$ we are looking for more detailed information  on
the growth of $T(r,f)$ and $n(r,c)$. The analog to Theorem
\ref{THM10} in connection with Theorems \ref{THM11} and
\ref{THM12} is as follows:

\begin{thm}\label{THM18}Suppose $f\in\Y_{\alpha,-1}$. Then
$f^\#(z)=O(|z|^{|\alpha|-1})$,
 \be{anzfktbeta0}n(r,\infty)\asymp\log r\quad{\rm and}\quad m(r,f)=O(\log
 r)\ee
hold; the same is true for every $c\in\C$ instead of $\infty$. In
particular we have
 \be{charfktbeta0}T(r,f)=N(r,f)+O(\log r)\asymp \log^2 r.\ee
\end{thm}

\proof For $\lambda>1$ we consider the annuli
$A_n:\lambda^{n-1}\le|z|<\lambda^n.$ By Theorem \ref{THM7}, each
$A_n$ contains at most $O((\lambda-1)^{-1})$ poles if $\lambda$ is
sufficiently close to $1$ (according to $\epsilon_0$ in Theorem
\ref{THM7}), and at least one pole if $\lambda$ is sufficiently
large (according to $\eta_0$ in Theorem \ref{THM7}). Thus
$n(r,\infty)\asymp \log r$ follows, and the same is true for any
value $c\in\C$ instead of $c=\infty$.\quad\Ende

\medskip\noindent{\sc Example.} Transcendental meromorphic solutions to algebraic
differential equations(\footnote{Yosida's contribution is known as
{\it Malmquist-Yosida Theorem} \cite{KY1,KY3}.}) $w'^n=R(z,w)$
have order of growth $\varrho\ge 1/3$ or else $\varrho=0$ (Bank
and Kaufman~\cite{BaKa2}, the author~\cite{NSt0}). An example for
the latter case is due to Bank and Kaufman~\cite{BaKa1} (slightly
modified): $\ds w'^2=\frac{4w(w^2-g_2/4)}{1-z^2}.$ One of its
solutions satisfies $f(\sin z)=\wp(z),$ where $\wp$ is the
Weierstrass P-function to the differential equation
$\wp'^2=4\wp(\wp^2-g_2/4)$, and $g_2>0$ is chosen in order that
$\wp$ has period lattice $\pi(\Z+i\Z)$. The zeros $\pm \cosh
\pi(k+\frac 12)$ and $\sqrt{g_2}/2$-points $\pm \cosh \pi k$ are
real, and the poles $\pm i\sinh(\pi k)$ and $-\sqrt{g_2}/2$-points
$\pm i\sinh\pi (k+\frac 12)$ are purely imaginary. From $\ds
f^\#(z)^2=\frac{4|f(z)||f(z)^2-g_2/4|}{|z^2-1|(1+|f(z)|^2)^2}$ and
the distribution of critical points follows $f^\#(z)=O(|z|^{-1})$,
and $f^\#(z)\asymp|z|^{-1}$ in $|\arg
z\pm\frac\pi4|<\frac\pi4-\epsilon$ and in $|\arg
z\pm\frac34\pi|<\frac\pi4-\epsilon$, hence $f\in\Y_{0,-1}$ by
Theorem \ref{THM3}. The $k$-th derivative of $f$ belongs to
$f\in\Y_{-k,-1}$. Any limit function satisfies $\f'(z)^2=\ds
-\frac{4\f(z)(\f(z)^2-g_2/4)}{(z+1)^2}$, hence has the form
$\f(z)=\wp(c+i\log(z+1))$ -- it is, of course, single-valued since
$\wp$ has period $\pi$, with essential singularity at $z=-1$.

\end{document}